\documentclass[10pt]{article}

\usepackage{a4wide}
\usepackage{amsmath}
\usepackage{amssymb}
\usepackage{amsthm}
\usepackage{amscd}
\usepackage{amsfonts}
\usepackage{graphicx}
\usepackage{latexsym}
\usepackage{eucal}
\usepackage{bm}
\usepackage{color}
\usepackage[multiple, stable]{footmisc}

\newcommand{\pp}[1]{{\bf(P#1)}}

\newcommand{\rr}[1]{{\bf Rule #1}}

\newenvironment{prf}[1][]{\par\noindent {\bf Proof#1:}\ }{\hfill$\blacksquare$
\par\vspace{11pt}}

\theoremstyle{plain} %\numberwithin{equation}{section}
\newtheorem{theorem}{Theorem}
\newtheorem{corollary}[theorem]{Corollary}
\newtheorem{conjecture}[theorem]{Conjecture}
\newtheorem{lemma}[theorem]{Lemma}
\newtheorem{proposition}[theorem]{Proposition}
\theoremstyle{definition}

\newtheorem{question}[theorem]{Question} \topmargin-2cm

\newtheorem{claim}[theorem]{Claim}

\newcommand{\Ccal}[0]{\ensuremath{{\mathcal K}}}

\newcommand{\Pcal}[0]{\ensuremath{{\mathcal P}}}
\newcommand{\Hcal}[0]{\ensuremath{{\mathcal H}}}
\newcommand{\Dcal}[0]{\ensuremath{{\mathcal D}}}
\newcommand{\Bcal}[0]{\ensuremath{{\mathcal B}}}

\newcommand{\Scal}[0]{\ensuremath{{\mathcal S}}}
\newcommand{\Tcal}[0]{\ensuremath{{\mathcal T}}}
\newcommand{\Ucal}[0]{\ensuremath{{\mathcal U}}}

\newcommand{\Wcal}[0]{\ensuremath{{\mathcal W}}}

\newcommand{\norm}[1]{\ensuremath{\|#1\|}}

\newcommand{\eR}[0]{\ensuremath{ \mathbb R}}

\newcommand{\eN}[0]{\ensuremath{ \mathbb N}}
\newcommand{\Zed}[0]{\ensuremath{ \mathbb Z}}

\newcommand{\Pee}[0]{\ensuremath{{\mathbb P}}}

\newcommand{\eps}[0]{\varepsilon}

\DeclareMathOperator{\diam}{diam}
\DeclareMathOperator{\area}{area}
\DeclareMathOperator{\vol}{vol}
\DeclareMathOperator{\Bi}{Bi}
\DeclareMathOperator{\cnr}{corner}
\DeclareMathOperator{\sde}{side}

\begin{document}

\title{Hamiltonicity of the random geometric graph}
\author{Michael Krivelevich\footnote{School of Mathematical Sciences, Tel Aviv University, Ramat Aviv 69978, Israel.
E-mail: {\tt krivelev@post.tau.ac.il, tobias@post.tau.ac.il}.}%
\footnote{Research supported in part by USA-Israel BSF grant 2006322, by
grant 1063/08 from the Israel Science Foundation, and by a Pazy memorial award.}
\and 
Tobias M\"uller\footnotemark[1]\footnote{Research partially supported through an ERC advanced grant.}
}

\maketitle

\begin{abstract}
Let $X_1,\dots, X_n$ be independent, uniformly random points from $[0,1]^2$. 
We prove that if we add edges between these points one by one by order of increasing edge length then, with probability tending to 1 
as the number of points $n$ tends to $\infty$, the resulting graph gets its first Hamilton cycle at exactly the same time it 
loses its last vertex of degree less than two.  
This answers an open question of Penrose and provides an analogue for the random geometric graph of a celebrated result 
of Ajtai, Koml{\'o}s and Szemer{\'e}di and independently of Bollob\'as 
on the usual random graph.
We are also able to deduce very precise information on the limiting probability that the random geometric graph is Hamiltonian
analogous to a result of Koml{\'o}s and Szemer{\'e}di on the usual random graph.
The proof generalizes to uniform random points on the $d$-dimensional hypercube 
where the edge-lengths are measured using the $l_p$-norm for some $1<p\leq\infty$.
The proof can also be adapted to show that, with probability tending to 1 
as the number of points $n$ tends to $\infty$, there are cycles of all lengths between $3$ and $n$ at the moment
the graph loses its last vertex of degree less than two.

\end{abstract}

\section{Introduction and statement of result}

Let $X_1,X_2,\dots \in [0,1]^2$ be a sequence of random points, chosen 
idenpendently and uniformly at random from $[0,1]^2$. 
For $n\in\eN$ and $r\geq 0$ the {\em random geometric graph} $G(n,r)$ has vertex set $V_n := \{X_1,\dots,X_n\}$ and an 
edge $X_iX_j \in E_n$ iff~$\norm{X_i-X_j} \leq r$.
The "hitting radius" $\rho_n(\Pcal)$ of an increasing graph property $\Pcal$ is the least $r$ such that
$G(n,r)$ satisfies $\Pcal$, i.e.:
\[
\rho_n(\Pcal) := \min\{ r \geq 0 : G(n,r) \text{ satisfies } \Pcal \}. 
\]
Recall that a graph is {\em Hamiltonian} if it has a Hamilton cycle (that is, a cycle that goes through all the vertices of the graph). 
An obvious necessary (but not sufficient) condition for the existence of a Hamilton cycle
is that the minimum degree is at least two.
In this paper we prove the following result:

\begin{theorem}\label{thm:main}
$\Pee\left[ \rho_n( \text{minimum degree }\geq 2 ) = \rho_n( \text{Hamiltonian} ) \right] \to 1$ as $n \to \infty$.
\end{theorem}

This answers a question of Penrose (see~\cite{penroseboek}, page 317) and provides an analogue for the random geometric graph of a celebrated result of Ajtai, Koml{\'o}s and Szemer{\'e}di~\cite{ajtaikomlosszemeredi85}
and independently of Bollob\'as~\cite{bollobas84} on the usual random graph.
Theorem~\ref{thm:main} can be stated alternatively as saying that if we add the edges between the points $X_1,\dots,X_n$
by order of increasing edge length then, with probability tending to 1 as $n\to\infty$, 
the resulting graph obtains its first Hamilton cycle at exactly the same time it loses its last vertex
of degree $<2$.
By combining Theorem~\ref{thm:main} with Theorem 8.4 from~\cite{penroseboek} 
(for completeness we have repeated the relevant special case of this theorem as Theorem~\ref{thm:penrosemindeg} below) we see that:

\begin{corollary}
Let $(r_n)_n$ be a sequence of nonnegative numbers, and write
$x_n := \pi n r^2_n - (\ln n + \ln \ln n)$.
Then:
\[
\lim_{n\to\infty} \Pee\left[ G(n,r_n) \text{ is Hamiltonian } \right]
= \left\{ \begin{array}{cl}
0 & \text{ if } x_n \to -\infty; \\
\exp\left[ -(\sqrt{\pi} + e^{-x/2})e^{-x/2}\right] & \text{ if } x_n \to x \in \eR; \\
1 & \text{ if } x_n \to +\infty.
 \end{array}\right..
\]
\label{cor:maincor}
\end{corollary}

Corollary~\ref{cor:maincor} provides an analogue for the random geometric graph of a result by
Koml{\'o}s and Szemer{\'e}di~\cite{komlosszemeredi83} on the limiting probability that the usual random graph
is Hamiltonian.

Previously, Petit~\cite{petitthesis} showed that if $(r_n)_n$ is chosen such that $r_n / \sqrt{\ln n / \pi n} \to \infty$ then
the random geometric graph $G(n,r_n)$ is Hamiltonian with probability tending to 1.
This was later sharpened by Diaz, Mitsche and P\'erez~\cite{diazmitscheperez07}
who showed that the same is true whenever $r_n \geq (1+\eps)\sqrt{\ln n / \pi n}$ with $\eps>0$ arbitrary (but fixed).
Our results are again an improvement and in a sense the final word on Hamiltonicity of the random geometric graph.
In Section~\ref{sec:concl} we shall nonetheless offer an idea for future research on Hamilton cycles in the
random geometric graph.

Since writing this paper it has come to our attention that both Balogh, Bollob\'as and Walters~\cite{jobal} and 
P\'erez and Wormald~\cite{perezcomm} have independently obtained essentially the same results at pretty much the same
time. Earlier Balogh, Kaul and Martin~\cite{martincomm} had proved Theorem~\ref{thm:main} in the case when
the Euclidean norm in the definition of the random geometric graph is replaced by the $l_\infty$-norm 
(i.e.~we add an edge between two points if their $l_\infty$-distance is less than $r$).

Our proof readily extends to arbitrary dimension and the $l_p$-norm for any $1<p\leq\infty$ 
(i.e.~the case where the points are i.i.d.~uniform on the
$d$-dimensional unit hypercube and $\norm{.}$ in the definition of the random geometric graph
is the $l_p$-norm), but we have chosen to focus on the two-dimensional random geometric graph
with the Euclidean norm for the sake of the clarity of our exposition.
In Section\ref{sec:generalise} we briefly explain the changes needed to make the proof work 
in the case of arbitrary dimension and the $l_p$-norm.

Our proof can also be adapted to show that, with probability tending to 1 as the number of points $n$ tends to infinity, 
the random geometric graph becomes {\em pancyclic} (i.e.~there are cycles of all lengths between 3 and $n$)
at precisely the same moment it first achieves minimum degree at least two.
In Section~\ref{sec:concl} we give a brief sketch the adaptations needed to squeeze this out of our proof.

Our proof of Theorem~\ref{thm:main} is inspired by the analysis in~\cite{diazmitscheperez07}.
Let us briefly outline the main steps in the proof.
%Penrose~\cite{penrosekconn} has already shown that $\rho_n(\text{minimum degree}\geq 2)
%= \rho_n(\text{2-connected})$ with probability tending to 1 as $n\to\infty$.
We pick an $r$ that is close to, but slightly less than $\rho_n(\text{minimum degree}\geq 2)$ and we dissect 
the unit square into squares of side $\eta r$ for a small constant $\eta$.
Next we consider an auxiliary graph $\Dcal$ consisting of the lower left hand corners of those squares
of our dissection that have at least 100 of the $X_i$s in them, where we connect two points of $\Dcal$ if their
distance is less than $r' := r(1-\eta\sqrt{2})$.
As it turns out, this auxiliary graph consists of one ``giant'' component
and a number of small components, that are cliques and are very far apart from each other.
Moreover, all of the $X_i$s are within distance $r$ of all the $\geq 100$ points in some
square of the auxiliary graph, except for a few clusters of ``bad'' points. These bad clusters
form cliques in the underlying random geometric graph, and these cliques are far apart.
We now construct a spanning tree $\Tcal$ of the giant component of $\Dcal$ that has maximum degree at most 26.
We increase $r$ to $\rho > r$ which is large enough for the random geometric graph to have minimum degree at least two, and we
construct the Hamilton cycle while performing a closed walk on $\Tcal$ that traveres every edge of $\Tcal$ exactly twice 
(once in each direction).
Each time the walk visits a node of $\Tcal$, the cycle visits a fresh $X_i$ inside the corresponding square. While doing this 
we are able to make small ``excursions'' to eat up the $X_i$s in squares belonging to non-giant components of $\Dcal$, 
the bad clusters and all the other $X_i$s.
% in squares belonging to non-giant components of $\Dcal$, in bad clusters and all other points.

\section{The proof}

Recall that a graph $G=(V,E)$ is $k$-connected if $|V| > k$ and $G\setminus S$ is connected for all sets $S \subseteq V$ of cardinality $|S| < k$.
Clearly, having minimum degree at least $k$ is a necessary condition for $k$-connectedness, and 2-connectedness is a necessary condition for Hamiltonicity.
In our proof of Theorem~\ref{thm:main} we shall rely on the following result of Penrose:

\begin{theorem}[\cite{penrosekconn}]\label{thm:penrosekconn}
For any (fixed) $k \in \eN$ it holds that:
\[ 
\Pee\left[ \rho_n( \text{minimum degree}\geq k) = 
\rho_n( \text{$k$-connected} ) \right] \to 1,
\] 
as $n\to\infty$.
\end{theorem}

Thanks to this last theorem, it suffices for us to show that
$\Pee[ \rho_n(\text{Hamiltonian}) = \rho_n(\text{2-connected})]\to 1$ in order to prove Theorem~\ref{thm:main}. 
We shall also make use of another result of Penrose. The following theorem is a reformulation of a special case of Theorem 8.4 from~\cite{penroseboek}.

\begin{theorem}\label{thm:penrosemindeg}
Let $(r_n)_n$ be a sequence of nonnegative numbers, and write
$x_n := \pi n r^2 - (\ln n + \ln \ln n)$.
Then:
\[
\lim_{n\to\infty} \Pee\left[ G(n,r_n) \text{ has minimum degree} \geq 2 \right]
= \left\{ \begin{array}{cl}
0 & \text{ if } x_n \to -\infty; \\
\exp\left[ - (\sqrt{\pi} + e^{-x/2})e^{-x/2}\right] & \text{ if } x_n \to x \in \eR; \\
1 & \text{ if } x_n \to +\infty.
 \end{array}\right.
\]
\end{theorem}

For $V \subseteq \eR^2$ and $r\geq 0$ we shall denote by $G(V,r)$ the (non-random) geometric graph with 
vertex set $V$ and an edge $vw\in E(G(V,r))$ iff~$\norm{v-w}\leq r$. 
The (non-random) geometric graphs $G(V,r)$ have been the subject of considerable research 
effort and they are often also called {\em unit disk graphs}.

For $0 < \eta < 1/\sqrt{2}$ and $r > 0$ let $\Hcal_\eta(r)$ denote the unit disk graph
$G(P_{\eta r},r')$ with vertex set $P_{\eta r} := [0,1]^2 \cap (\eta r)\Zed^2$ (that is,
$P_{\eta r}$ is the set of all points in $[0,1]^2$ whose coordinates are integer 
multiples of $\eta r$) and threshold distance $r' := r(1-\eta\sqrt{2})$.

Now suppose that we are also given an arbitrary set $V\subseteq [0,1]^2$ of points.
We shall call a vertex $p\in \Hcal_\eta(r)$ {\em dense} with respect to $V$ if the square $p + [0,\eta r)^2$ contains at least
$100$ points of $V$.
If a vertex is not dense we will call it {\em sparse}.
If  all neighbours of $p$ in $\Hcal_\eta(r)$ are sparse (i.e.~if $q$ is sparse for all $q\in B(p, r')\cap\Hcal_\eta(r)$) then we shall say that
$p$ is {\em bad}.
 
Let $\Dcal_\eta(V,r)$ denote the subgraph of $\Hcal_\eta(r)$ induced by the
dense points, and let $\Bcal_\eta(V,r)$ denote the subgraph induced by the 
bad points. 
Part of the proof of Theorem~\ref{thm:main} will be to show that
if $V=\{X_1,\dots,X_n\}$ and 
$r$ is chosen close to, but slightly smaller than, $\rho_n(\text{2-connected})$ then
$\Hcal_\eta(r), \Dcal_\eta(V,r)$ and $\Bcal_\eta(V,r)$ have a number of 
desirable properties (with probability tending to 1).
This will then allow us to finish the proof of our main theorem by purely deterministic arguments.  
Here is a list of these desirable properties
(here and throughout the rest of the paper ``component'' will always mean a connected component, and
diameter will always refer to the geometric diameter of a point set as opposed to the graph diameter):

\begin{enumerate}
\item[\pp{1}] If $\Ccal$ is a component of $\Dcal_\eta(V,r)$ then it
either has (geometric) diameter $\diam(\Ccal) < r'$ or
$\diam(\Ccal) > 1000r$;
\item[\pp{2}] If $\Ccal_1,\Ccal_2$ are two distinct components of $\Dcal_\eta(V,r)$ with
(geometric) diameters $\diam(\Ccal_1),\diam(\Ccal_2) < r'$ and 
$p_1\in \Ccal_1, p_2\in\Ccal_2$ then $\norm{p_1-p_2} > 1000r$;
\item[\pp{3}] If $p_1\in\Ccal$ for some component $\Ccal$ of $\Dcal_\eta(V,r)$ with $\diam(\Ccal) < r'$ and $p_2\in\Bcal_\eta(V,r)$ is 
bad then $\norm{p_1-p_2} > 1000 r$;
\item[\pp{4}] If $p,q \in \Bcal_\eta(V,r)$ are bad, then either
$\norm{p-q} < r'$ or $\norm{p-q} > 1000 r$;
\item[\pp{5}] If $p_1,p_2 \in \Dcal_\eta(V,r)$ and 
$\norm{p_1-p_2} < 25r$ and neither of $p_1$ or $p_2$ lies in a component
of (geometric) diameter $< r'$ then there is a 
$p_1p_2$-path in $\Dcal_\eta(V,r)$ that stays inside $B(p_1,100r)$;
\item[\pp{6}] $\Dcal_\eta(V,r)$ has exactly one component $\Ccal$ of (geometric) diameter $\diam(\Ccal) \geq r'$.
\end{enumerate}

We will say that a sequence of events $(A_n)_n$ holds {\em with high probability} (w.h.p.) if
$\Pee( A_n ) \to 1$ as $n\to\infty$.
The following proposition takes care of the probabilistic part of the proof of Theorem~\ref{thm:main}:

\begin{proposition}\label{prop:Hstruct}
Set $r_n := \sqrt{\ln n / \pi n}$ and $V_n := \{X_1,\dots,X_n\}$. 
If $\eta > 0$ is sufficiently small (but fixed) then $\Hcal_\eta(r_n), \Dcal_\eta(V_n,r_n)$ and $\Bcal_\eta(V_n,r_n)$
satisfy properties \pp{1}-\pp{6} w.h.p.
\end{proposition}

Together with the following deterministic result and the two mentioned results by Penrose, Proposition~\ref{prop:Hstruct}
gives  Theorem~\ref{thm:main}. 

\begin{theorem}\label{thm:deterministic} 
Suppose that $0 < \eta < 1/\sqrt{2}, V \subseteq \eR^2$ and $r>0$ are such that $\Dcal_\eta(V,r)$ and $\Bcal_\eta(V,r)$ satisfy \pp{1}-\pp{6}, and that $r \leq \rho \leq 2r$ is 
such that $G(V,\rho)$ is 2-connected.
Then $G(V,\rho)$ is also Hamiltonian.
\end{theorem}

We postpone the proofs of Proposition~\ref{prop:Hstruct} and Theorem~\ref{thm:deterministic} and we first briefly explain how 
they imply Theorem~\ref{thm:main}.

\vspace{11pt}

\begin{prf}[ of Theorem~\ref{thm:main}] 
Let us write $r_n := \sqrt{\ln n / \pi n}$ and $\sigma_n :=  \rho_n(\text{2-connected})$.
By Theorem~\ref{thm:penrosekconn} and the fact that 2-connectedness is a necessary condition
for Hamiltonicity it suffices to show that $G(n,\sigma_n)$ is Hamiltonian with high probability. 
Theorem~\ref{thm:penrosekconn} together with Theorem~\ref{thm:penrosemindeg}
show that, with high probability, $G(n,r_n)$ is not 2-connected and $G(n,2r_n)$ is 2-connected.
In other words, $r_n < \sigma_n \leq 2r_n$ with high probability.
By Proposition~\ref{prop:Hstruct} we can fix an $\eta \in (0,1/\sqrt{2})$ such that 
properties \pp{1}-\pp{6} hold for 
$\Hcal_\eta(r_n), \Dcal_\eta(V_n,r_n), \Bcal_\eta(V_n, r_n)$ with high probability, where $V_n := \{X_1,\dots,X_n\}$.
Thus, with high probability, Theorem~\ref{thm:deterministic} applies to $\eta, V = \{X_1,\dots,X_n\}, r=r_n, \rho=\sigma_n$, and 
$G(n,\sigma_n)$ is indeed Hamiltonian with high probability.
\end{prf}

Our next step is to prove Proposition~\ref{prop:Hstruct}. 
We will say that a point or set is {\em within $s$ of the sides} (of $[0,1]^2$)
if it is contained in
\[ \sde(s) := \{ z\in[0,1]^2 : z_x\in [0,s)\cap (1-s,1] \text{ or }
z_y \in [0,s)\cap (1-s,1] \},
\]
and we will say it is {\em within $s$ of the corners} (of $[0,1]^2$) if
it is contained in
\[ 
\cnr(s) := \{ z\in[0,1]^2 : z_x\in [0,s)\cap (1-s,1] \text{ and }
z_y \in [0,s)\cap (1-s,1] \}.
\]
The following lemma provides an observation 
that is pivotal in the proof of Proposition~\ref{prop:Hstruct}.

\begin{lemma}\label{lem:Sdense}
Set $r_n := \sqrt{\ln n / \pi n}$ and $V_n := \{X_1,\dots,X_n\}$.
For every $\eps > 0$ there exists an $\eta_0 = \eta_0(\eps) > 0$ such that 
for any fixed $0<\eta<\eta_0$ the following statements hold w.h.p.:
\begin{enumerate}
\item\label{itm:Sdense1} For every $\Scal \subseteq \Hcal_\eta(r_n)$ with 
$|\Scal| > (1+\eps)\pi\eta^{-2}$ and $\diam(\Scal)<10^5r_n$,
there exists a $q\in\Scal$ that is dense wrt.~$V_n$;
\item\label{itm:Sdense2} For every $\Scal \subseteq \Hcal_\eta(r_n)\cap\sde(10^5r_n)$ with 
$|\Scal| > (1+\eps)\frac{\pi}{2}\eta^{-2}$ and $\diam(\Scal) < 10^5r_n$, there exists 
a $q\in\Scal$ that is dense wrt.~$V_n$;
\item\label{itm:Sdense3} For every $\Scal \subseteq \Hcal_\eta(r_n)\cap\cnr(10^5r_n)$ with  
$|\Scal| > \eps\eta^{-2}$ and $\diam(\Scal) < 10^5r_n$, there
exists a $q\in\Scal$ that is dense wrt.~$V_n$.
\end{enumerate}
\end{lemma}

In the proof of Lemma~\ref{lem:Sdense} we shall make use of the following incarnation of the Chernoff-Hoeffding bound.
A proof can for instance be found in~\cite{penroseboek}, on page 16.

\begin{lemma} Let $Z$ be a $\Bi(n,p)$-distributed random variable, and $k\leq \mu:= np$. Then
\[ 
\Pee( Z \leq k ) 
\leq 
\exp\left[ -\mu H(k/\mu) \right],
\]
where $H(x) := x\ln x - x + 1$.
\label{lem:chernoff} 
\end{lemma}

\begin{prf}[ of Lemma~\ref{lem:Sdense}] 
Let us choose $\eta_0 := \eps / 10^{6}$ and fix an arbitrary $0<\eta<\eta_0$. 
Our choice of $\eta_0$ guarantees that $4\lceil 2\cdot 10^5/\eta\rceil < \eps \eta^{-2} / 2$
(we can assume w.l.o.g.~that $\eps < 1$).

Let $\Ucal$ denote the collection of all
$\Scal \subseteq \Hcal_\eta(r_n)$ that satisfy the conditions for part~\ref{itm:Sdense1}.
Let us first count the number of sets in $\Scal$.
To this end, observe that if $p\in\Scal$ and $\diam(\Scal)<10^5r_n$ then
$\Scal \subseteq p+(-10^5r_n,10^5r_n)^2$.
Notice that
\[ 
\left|\Hcal_\eta(r_n) \cap \left(p+(-10^5r_n,10^5r_n)^2\right)\right| \leq \lceil 2\cdot 10^5/\eta\rceil^2, 
\]
for any $p\in\eR^2$.
This shows that if $N(p)$ denotes the number of $\Scal\in\Ucal$ that contain $p$, then 
\begin{equation}\label{eq:Nycnt} 
N(p) \leq 2^{\lceil 2\cdot 10^5/\eta\rceil^2}, 
\end{equation}
and, since this constant upper bound on $N(p)$ holds for all $p\in\Hcal_\eta(r_n)$, it follows that 
\begin{equation}\label{eq:Scalcnt} 
|\Ucal| \leq \sum_{p\in\Hcal_\eta(r_n)} 2^{N(p)} = O( |\Hcal_\eta(r_n)| ) = O( r_n^{-2} ) = O( n / \ln n ). 
\end{equation}
Now pick an arbitrary $\Scal\in\Ucal$, and let
$\Scal'\subseteq \Scal$ be the set of those $q\in\Scal$ for which
$q+[0,\eta r_n)^2 \subseteq [0,1]^2$.
Since $\Scal \subseteq p+(-10^5r_n,10^5r_n)^2$ for any $p\in\Scal$, we have that 
$|\Scal \setminus \Scal'| \leq 4 \cdot \lceil 2\cdot 10^5/\eta\rceil$.
And hence, by choice of $\eta_0$, we have 
\begin{equation}\label{eq:Scalcard}
|\Scal'| > (1+\frac{\eps}{2})\pi\eta^{-2}. 
\end{equation}
Let $Z := |\{X_1,\dots, X_n\} \cap (\bigcup_{q\in\Scal'} q+[0,\eta r)^2 )|$ denote the (random) number of 
$X_i$ that fall into one of the squares $q+[0,\eta r_n)^2$ with $q\in\Scal'$.
Then $Z \sim \Bi( n, |\Scal'|\eta^2r_n^2 )$.
Appealing to Lemma~\ref{lem:chernoff}:
\begin{equation}\label{eq:sparseforall}
\begin{split}
\Pee\left[ q \text{ is sparse for all } q\in \Scal \right] 
& \leq 
\Pee\left[ q \text{ is sparse for all } q\in \Scal' \right] \\
& \leq
\Pee\left[ Z \leq 99|\Scal'| \right] \\
& \leq 
\exp[ - n|\Scal'|\eta^2r^2_n H(99 / \eta^2 nr^2_n) ],
\end{split}
\end{equation}
where $H(x) = x\ln x - x + 1$.
Now notice that, by~\eqref{eq:Scalcard}
\begin{equation}\label{eq:muexpr}
 n|\Scal'|\eta^2r^2_n H(99 / \eta^2 nr^2_n)
> (1+\frac{\eps}{2})\pi n r_n^2 H(99 / \eta^2 nr^2_n)
= (1+\frac{\eps}{2}+o(1)) \ln n,
\end{equation}
where the last equality holds by the choice of $r_n = \sqrt{\ln n / \pi n}$ and the fact
that $H(x)\to 1$ as $x\downarrow 0$ (note $99 / (\eta^2 n r^2_n) = O( 1 / \ln n ) \to 0$). 
Combining~\eqref{eq:Scalcnt},~\eqref{eq:sparseforall} and~\eqref{eq:muexpr}, the union bound now gives us that
\[ 
\begin{split}
\Pee\left[
\exists \Scal \in \Ucal \text{ such that }
q \text{ sparse for all } q \in\Scal \right]
& \leq
\sum_{\Scal\in\Ucal}
\Pee\left[
q \text{ sparse for all } q \in \Scal \right] \\
& \leq 
|\Ucal| n^{-1-\frac{\eps}{2}+o(1)} \\
& = 
o(1), 
\end{split}
\]
which proves part~\ref{itm:Sdense1}.

Now let $\Ucal_{\sde} \subseteq \Ucal$ denote the collection of
all $\Scal\subseteq\Hcal_\eta(r_n)$ that satisfy the conditions of
part~\ref{itm:Sdense2} of the lemma.
Noticing that
\[ 
|\Hcal_\eta(r_n) \cap \sde( 10^5r_n )| \leq 4 \cdot \lceil 10^5 / \eta\rceil \cdot (1 + 1 / \eta r_n)
= O( 1/ r_n ) = O( \sqrt{n / \ln n } ), 
\]
and reusing~\eqref{eq:Nycnt}, we
see that:
\begin{equation}\label{eq:Tcalcnt}
|\Ucal_{\sde}| \leq \sum_{p\in\Hcal_\eta(r)\cap \sde(10^5r_n)} 2^{N(p)} = O( \sqrt{n / \ln n} ).
\end{equation} 
Now let $\Scal\in\Ucal_{\sde}$ be arbitrary, and let $\Scal'\subseteq\Scal$ be those
$q\in\Scal$ for which $q+[0,\eta r_n)^2 \subseteq [0,1]^2$.
Again we have $|\Scal\setminus\Scal'| \leq 4\lceil 2\cdot 10^5 \eta^{-1}\rceil < \frac{\eps}{2}\eta^{-2}$, so that this time
$|\Scal'| > (1+\frac{\eps}{2})\frac{\pi}{2}\eta^{-2}$.
The inequality~\eqref{eq:sparseforall} is still valid and, analogously to~\eqref{eq:muexpr}, we 
now have $n|\Scal'|\eta^2r^2_n H(99 / \eta^2 nr^2_n)
> (\frac12+\frac{\eps}{4}+o(1)) \ln n$.
Combining these observations with~\eqref{eq:Tcalcnt}, the union bound thus gives:
\[
\Pee\left[ \exists \Scal\in\Ucal_{\sde}\text{ such that } q \text{ sparse for all } q\in\Scal\right]
\leq 
|\Ucal_{\sde}| n^{-\frac12-\frac{\eps}{4}+o(1)} = o(1),
\]
proving part~\ref{itm:Sdense2} of the lemma.

Finally, let $\Ucal_{\cnr}$ denote the collection of sets $\Scal\subseteq \Hcal_\eta(r_n)$ that satisfy
the conditions of part~\ref{itm:Sdense3} of the lemma.
Notice that $|\Hcal_\eta(r_n)\cap \cnr(10^5r_n)| \leq 4\lceil 10^5\eta^{-1}\rceil^2 = O(1)$.
Therefore, also:
\begin{equation}\label{eq:Ucalcard}
 |\Ucal_{\cnr}| \leq \sum_{p\in\Hcal_\eta(r_n)\cap \cnr(10^5r_n)} 2^{N(p)} = O(1). 
\end{equation}
Pick an arbitrary $\Scal\in\Ucal_{\cnr}$ and let $\Scal'\subseteq\Scal$ be the set of those
$q\in\Scal$ for which $q+[0,\eta r_n)^2 \subseteq [0,1]^2$.
Then $|\Scal'| > \frac{\eps}{2}\eta^{-2}$. 
Again the inequality~\eqref{eq:sparseforall} is still valid and, analogously to~\eqref{eq:muexpr}, we 
now have $n|\Scal'|\eta^2r^2_n H(99 / \eta^2 nr^2_n)
> (\frac{\eps}{2}+o(1)) \ln n$.
Combining this with~\eqref{eq:Ucalcard}, the union bound gives:
\[ 
\Pee\left[ \exists \Scal\in\Ucal_{\cnr} \text{ such that } q \text{ sparse for all }
q \in \Scal \right] 
\leq |\Ucal_{\cnr}| n^{-\frac{\eps}{2}+o(1)} = o(1).
\]
This proves part~\ref{itm:Sdense3} of the lemma.
\end{prf}

We say that a set $A \subseteq \eR^2$ is a {\em Boolean combination} of the sets $A_1,\dots, A_n \subseteq \eR^2$ if  
$A$ can be constructed from $A_1,\dots,A_n$ by means of any number of compositions of the operations intersection, union and complement. 
Recall that a {\em halfplane} is a set of the form $H(a,b) := \{ z \in \eR^2 : z\cdot a \leq b\}$ for some
vector $a \in \eR^2 \setminus\{0\}$ and constant $b\in\eR$.

\begin{lemma} There exists a constant $C>0$ such that the following holds for 
all $\eta, r > 0$.
For every $A \subseteq\eR^2$ with $\diam(A) < 10^5r_n$ that is a Boolean combination of at most $1000$ halfplanes and balls of radius $\leq r$, we have that $|A\cap\Hcal_\eta(r)| \geq \area(A\cap [0,1]^2) / (\eta r)^2 - C \eta^{-1}$.
\label{lem:Sarea}
\end{lemma}

\begin{prf} Set $C := 10^9$, and let $\eta, r>0$ be arbitrary.
Let $A \subseteq \eR^2$ be an abitrary set that satisfies the two conditions from the lemma.
Let 
\[ 
A' := \{ z \in \eR^2 : B(z;\eta r\sqrt{2}) \subseteq A\cap [0,1]^2 \}.
\] 
In other words, $A'\subseteq A\cap [0,1]^2$ is the set of all $z$ that are distance at least $\eta r\sqrt{2}$ away from the boundary of $A\cap [0,1]^2$.
Observe that if $q + [0,\eta r)^2$ intersects $A'$ then it is completely contained in $A$.
Because the squares $q + [0,\eta r)^2 : q \in \Hcal_\eta(r)$ are disjoint and cover $[0,1]^2$ this shows that
\[
|A\cap\Hcal_\eta(r)| \geq \area(A') / (\eta r)^2.
\]
It thus suffices to bound the area of $A\cap [0,1]^2\setminus A'$.
The set $A\cap[0,1]^2$ is also a Boolean combination of halfplanes and balls of radius $\leq r$, this time
at most 1004 of them.
Let $H_1 = H(a_1,b_1), \dots, H_m = H(a_m,b_m)$
and $B_1 := B(z_1;s_1), \dots, B_k = B(z_k;s_k)$ 
denote the halfplanes and disks that $A\cap [0,1]^2$ is a Boolean combination of, where
$m+k \leq 1004$ and $s_1,\dots, s_k \leq r$. 
We can assume w.l.o.g.~that $\norm{a_i} = 1$ for $i=1,\dots,m$.
Pick an arbitrary $z_0\in A$. Because $\diam(A) < 10^5r$ we have $A\subseteq B(z_0,10^5r)$.
Let us now observe that if $z\in (A\cap[0,1]^2)\setminus A'$ then
$z$ lies within distance $\eta r\sqrt{2}$ of the boundary 
of one of the sets $H_i$ or one of the $B_j$.
This implies that
\[ 
A\cap [0,1]^2 \setminus A'
\subseteq H_1'\cup \dots  \cup H_m' \cup B_1'\cup \dots \cup B_k', \]
where $H_i' := B(z_0,10^5r) \cap H(a_i,b_i+\eta r\sqrt{2})\setminus H(a_i, b_i-\eta r\sqrt{2})$ for 
$i=1,\dots,m$ and 
$B_j' := B( z_j; s_j + \eta r\sqrt{2} ) \setminus B( z_j; s_j - \eta r\sqrt{2} )$ for $j=1,\dots,k$.
Now notice that $\area( H_i') \leq (2\cdot 10^5r)\times(2\eta r\sqrt{2}) = 4 \cdot 10^5 \eta r^2\sqrt{2}$ and 
$\area( B_i') = \pi ( ( s_i + \eta r\sqrt{2} )^2 - (s_i - \eta r\sqrt{2})^2 ) =
4 \pi s_i \eta r\sqrt{2} \leq 4 \pi \eta r^2 \sqrt{2}$.
Thus
\[ 
\area( A\cap[0,1]^2\setminus A' ) \leq 1004 \cdot 4\cdot 10^5 \eta r^2\sqrt{2} \leq C \eta r^2, 
\]
which gives $|A\cap \Hcal_\eta(r)| \geq \left( \area(A\cap [0,1]^2) - C \eta r^2 \right) / (\eta r)^2 
= \area(A\cap [0,1]^2) / (\eta r)^2 - C \eta^{-1}$ as required.
\end{prf}

\begin{prf}[ of Proposition~\ref{prop:Hstruct}]
Set $V_n=\{X_1,\dots,X_n\}, r_n=\sqrt{\ln n / \pi n}$.
We will show how
Lemma~\ref{lem:Sdense} can be applied to show that each of 
the statements \pp{1}-\pp{6} hold with high probability for $\Hcal_\eta(r_n), \Dcal_\eta(V_n,r_n), \Bcal_\eta(V_n,r_n)$ 
if $\eta$ is chosen sufficiently small.

Set $\eps := 1/1000$. Fix an $0<\eta< \eta_0( \eps )$, where
$\eta_0$ is as in Lemma~\ref{lem:Sdense}, that is also small enough for the following three 
inequalities to hold:
\begin{equation}
\begin{split}
(1+\frac{1}{100})(1-\eta\sqrt{2})^2\pi - C \eta & > (1+\eps)\pi, \\
(1+\frac{1}{100})(1-\eta\sqrt{2})^2\frac{\pi}{2} - C\eta & > (1+\eps)\frac{\pi}{2}, \\
(1-\eta\sqrt{2})^2\frac{\pi}{4} - C\eta & > \eps,
\end{split}
\label{eq:etadef}
\end{equation}
where $C$ is the constant from Lemma~\ref{lem:Sarea}.

For any $r>0$ (and $\eta,\eps$ as chosen above) let 
$\Ucal(r)$ denote the set of all $\Scal\subseteq\Hcal_\eta(r)$ for 
which $\diam(\Scal) < 10^5r$ and either $|\Scal| > (1+\eps)\pi\eta^{-2}$, or
$\Scal\subseteq\sde(10^5r)$ and $|\Scal|>(1+\eps)\frac{\pi}{2}\eta^{-2}$, or
$\Scal\subseteq\cnr(10^5r)$ and $|\Scal|>\eps\eta^{-2}$. 
By Lemma~\ref{lem:Sdense} it holds with high probability that any $\Scal\in\Ucal(r_n)$ contains a point that is dense wrt.~$V_n$.
To prove the proposition it thus suffices to show that for any $V\subseteq [0,1]^2$ and $0<r<10^{-10}$ that are such that
each $\Scal\in\Ucal(r)$ contains a point that is dense wrt.~$V$
%these three statements hold for $\Hcal_\eta(r)$ wrt.~$V$ (where $\eps,\eta$ are as chosen above), 
the properties \pp{1}-\pp{6} hold for $\Hcal_\eta(r), \Dcal_\eta(V,r)$ and $\Bcal_\eta(V,r)$ 
(with $\eta,\eps$ as chosen above). 
Let us thus pick such a $V\subseteq[0,1]^2$ and $0<r<10^{-10}$ for which 
every $\Scal\in\Ucal(r)$ contains a point that is dense wrt.~$V$.

\noindent {\bf Proof that \pp{1} holds:}
Aiming for a contradiction, suppose there is some component $\Ccal$ of $\Dcal_\eta(n,r)$ with (geometric) diameter
$r' \leq \diam(\Ccal) \leq 1000r$.
Let us pick points $p_L,p_R,p_T,p_B\in\Ccal$, where 
$p_L$ is a point of $\Ccal$ with smallest $x$-coordinate amongst all points of $\Ccal$, $p_R$ is a point of
$\Ccal$ with biggest $x$-coordinate, $p_B$ is a point of
$\Ccal$ with smallest $y$-coordinate and $p_T$ is a point of $\Ccal$ 
with biggest $y$-coordinate (note these points need not be distinct or unique). See Figure~\ref{fig:smallcomps}
for an illustration.
Since $\diam(\Ccal) \geq r'$, we have either
$(p_R)_x-(p_L)_x \geq r'/\sqrt{2}$ or 
$(p_T)_y-(p_B)_y \geq r'/\sqrt{2}$.
Without loss of generality, let us assume that $(p_R)_x-(p_L)_x \geq r'/\sqrt{2}$.
For $z\in\eR^2$ and $s\geq 0$ let us set:

\[ \begin{array}{l}
B_L(z,s) := \{ z' \in B(z, s) : z_x' < z_x \}, \quad
B_R(z,s) := \{ z' \in B(z, w) : z_x' > z_x \}, \\
B_B(z,s) := \{ z' \in B(z, r') : z_y' < z_y \}, \quad
B_T(z,s) := \{ z' \in B(z, r') : z_y' > z_y \}.
\end{array} \]
Now define
\[ A := B_L(p_L;r')\cup B_R(p_R,r')\cup B_B(p_B,r')\cup B_T(p_T,r'), \]
and let $\Scal := A \cap \Hcal_\eta(r)$ denote the set 
of all points of $\Hcal_\eta(r)$ that fall inside $A$.
Let us observe that, since $\Ccal$ is a component of $\Dcal_\eta(n,r)$, the
set $\Scal$ cannot contain any dense $q$.
We also note that $\diam(A) < 10^5r$ and $A$ is a Boolean combination of $\leq 1000$ halfspaces
and balls of radius $\leq r$.

Let us define
\[ 
B_B' := \{ z \in B_B(p_B,r') : (p_L)_x < z_x < (p_R)_x \}, \quad
B_T' := \{ z \in B_T(p_T,r') : (p_L)_x < z_x < (p_R)_x \}. 
\]
Then the sets $B_L(p_L,r'), B_R(p_R,r'), B_B', B_T'$ are disjoint (see Figure~\ref{fig:smallcomps}).
\begin{figure}[!hbt]
\begin{center}
\input{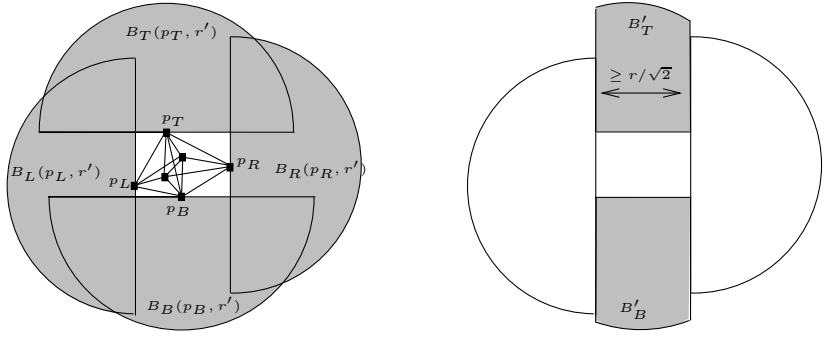}
\end{center}
\caption{$A$ has area at least $(1+\frac{1}{2\sqrt{2}})\pi(r')^2$.}
\label{fig:smallcomps}
\end{figure}
Now observe that, because $(p_R)_x - (p_L)_x > r'/\sqrt{2}$ and 
$(p_L)_x \leq (p_T)_x \leq (p_R)_x$, 
the area of $B_T'$ is at least a fraction $(r'/\sqrt{2}) / (2r') = 1/2\sqrt{2}$ of the area of $B_T(p_T,r')$.
Similarly $B_B'$ has at least $1/2\sqrt{2}$ of the area of $B_B(p_B,r')$.
In other words, $\area( B_B'), \area( B_T' ) \geq \frac{1}{4\sqrt{2}}\pi (r')^2$,
and thus
\begin{equation}\label{eq:H3Hlb}
\area(A) \geq (1+\frac{1}{2\sqrt{2}}) \pi (r')^2.
\end{equation}
First suppose that $A$ is completely contained in $[0,1]^2$.
In this case, Lemma~\ref{lem:Sarea} tells us that
\[ \begin{split}
 |\Scal| 
& \geq \area(A) / (\eta r)^2 - C \eta^{-1} \\
& \geq (1+\frac{1}{2\sqrt{2}})(1-\eta\sqrt{2})^2\pi \eta^{-2} - C\eta^{-1} \\
& > (1+\eps)\pi\eta^{-2},
\end{split}
\]
where the last inequality holds by~\eqref{eq:etadef}. We see that $\Scal \in \Ucal(r)$.
But then there must be a dense $q\in\Scal$, which contradicts that $\Ccal$ is a component of $\Dcal_\eta(V,r)$.

Now assume that one of the points $p_L,p_R,p_B,p_T$ is within distance
$r'$ of one of the sides of $[0,1]^2$, but none of these points is an element of $\cnr(r')$.
Then certainly $\Scal \subseteq \sde(10^5r)$.
Also note that at least one of $B_B', B_T'$ must be completely contained in $[0,1]^2$. 
Moreover, at least half the area of $B_L(p_L,r')\cup B_R(p_R,r')$ lies in $[0,1]^2$.
(If the points are close to $\{0\}\times[0,1]$ then
$B_R(p_R,r') \subseteq [0,1]^2$. If they are close to $\{1\}\times[0,1]$ then 
$B_L(p_L,r')\subseteq [0,1]^2$.
If they are close to $[0,1]\times\{0\}$ then the top halves of $B_L(p_L,r')$ and $B_R(p_R,r')$ lie completely in $[0,1]^2$. 
If they are close to $[0,1]\times\{1\}$ then the bottom halves of $B_L(p_L,r')$ and
$B_R(p_R,r')$ are completely contained in $[0,1]^2$.) 
Hence 
\[ \area( A \cap [0,1]^2 ) \geq (1 + \frac{1}{2\sqrt{2}})\frac{\pi}{2}(r')^2. \]
Using Lemma~\ref{lem:Sarea} and~\eqref{eq:etadef} we find:
\[ \begin{split}
|\Scal| 
& \geq \area(A \cap [0,1]^2 ) / (\eta r)^2 - C\eta^{-1} \\
& \geq (1 + \frac{1}{2\sqrt{2}})(1-\eta\sqrt{2})^2\frac{\pi}{2}\eta^{-2} - C\eta^{-1} \\
& > (1+\eps)\frac{\pi}{2}\eta^{-2}.
\end{split}
\]
Again we see that $\Scal\in\Ucal(r)$. So again at least one $q\in \Scal$ must be dense, which again contradicts that $\Ccal$ was a component
of $\Dcal_\eta(V,r)$.

Finally assume that one of the 4 points is an element of $\cnr(r')$.
Clearly $\Scal\subseteq \cnr(10^5r)$.
Also note that at least one of $B_B', B_T'$ is completely contained in $A\cap[0,1]^2$.
Hence, by Lemma~\ref{lem:Sarea} and~\eqref{eq:etadef}:
\[ \begin{split}
|\Scal| 
& \geq \area(A \cap [0,1]^2 ) / (\eta r)^2 - C\eta^{-1} \\
& \geq (1-\eta\sqrt{2})^2\frac{1}{4\sqrt{2}}\eta^{-2} - C\eta^{-1} \\
& > \eps\eta^{-2}.
\end{split}
\]
And again this implies $\Scal\in\Ucal(r)$, which in turn implies the existence of a dense $q\in\Scal$, which cannot be.
We can thus conclude that no component $\Ccal$ of diameter $r' \leq \diam(\Ccal) \leq 1000 r$ exists
in $\Dcal_\eta(V,r)$.

\noindent {\bf Proof that \pp{2} holds:}
Suppose that $\Ccal_1,\Ccal_2$ are distinct components of $\Dcal_\eta(V,r)$, both with
diameter $\leq r'$, such that there exists a point in $\Ccal_1$ and a point in $\Ccal_2$ 
with distance at most $1000r$ between them.
Now set $\Ccal := \Ccal_1\cup\Ccal_2$. Then $r'<\diam(\Ccal) \leq 1002r$
(to see the lower bound, note that any point in $\Ccal_1$ has distance $> r'$ to any point of $\Ccal_2$ as they are in distinct components).
Let $p_L \in \Ccal$ be a point of smallest $x$-coordinate, let $p_R\in\Ccal$ be a point of largest $x$-coordinate,
let $p_B\in\Ccal$ be a point of smallest $y$-coordinate, let $p_T$ be a point of largest $y$-coordinate and set
$A := B_L(p_L;r')\cup B_R(p_R,r')\cup B_B(p_B,r')\cup B_T(p_T,r')$.
Then $\Scal := A\cap \Hcal_\eta(r)$ cannot contain any dense point
(if, for example, $B_L(p_L,r')$ were to contain a dense point, then this point would lie in the same component
of $\Dcal_\eta(V,r)$ as $p_L$ and have a smaller $x$-coordinate than $p_L$).
We can now proceed as in the proof of \pp{1} to arrive at a contradiction.

\noindent {\bf Proof that \pp{3} holds:} Suppose that $\Ccal_1$ is a component of $\Dcal_\eta(V,r)$ with $\diam(\Ccal_1) < r'$
and that $p\in\Bcal_\eta(V,r)$ is a bad point that is at distance $< 1000r$ to some point of $\Ccal_1$.
Let us set $\Ccal := \Ccal_1\cup\{p\}$.
Then $r' < \diam(\Ccal) < 1001r$ (to see the lower bound, note that 
any bad point has distance $>r'$ to any dense point).
Defining $p_L,p_R,p_B,p_T, A$ and $\Scal$ as in the proofs of \pp{1} and \pp{2}, we see that $\Scal$ again cannot
contain any dense point. We again arrive at a contradiction by proceeding as in the proof of \pp{1}.

\noindent {\bf Proof that \pp{4} holds:} Suppose that $p_1,p_2\in\Bcal_\eta(V,r)$ are bad and that $r' \leq \norm{p_1-p_2} \leq 1000r$.
Setting $\Ccal := \{p_1,p_2\}$ and repeating the same argument again gives a contradiction. 

\noindent {\bf Proof that \pp{5} holds:}
Suppose there exist $p_1,p_2\in \Dcal_\eta(n,r)$ with $\norm{p_1-p_2} < 25r$ such that both points are in components
of diameter $\geq r'$, and there is no $p_1p_2$-path that stays inside $B(p_1,100r)$.
By \pp{1} $p_1,p_2$ must each be in a component of diameter $> 1000 r$.
For $k=25,\dots,70$ let $S_k := p_1 + [-kr,kr]^2$ denote a square of side length $2kr$ with center $p_1$, and 
let $R_k := S_{k} \setminus S_{k-1}$ for $k=26,\dots,70$.
Consider the subgraph $\tilde{\Dcal}$ of $\Dcal_\eta(V,r)$ induced
by the points of $\Dcal_\eta(V,r)$ that lie inside $S_{70}$.
Observe that $p_1$ and $p_2$ must lie in distinct components
$\tilde{\Ccal_1}, \tilde{\Ccal_2}$ of $\tilde{\Dcal}$ (otherwise there is a $p_1p_2$-path that stays inside
$S_{70} \subseteq B(p_1,100r)$) and that $R_k$ must contain a point of $\tilde{\Ccal_1}$ and of $\tilde\Ccal_2$ for each 
$k=26, \dots, 70$ (otherwise, if $\tilde\Ccal_j$ misses $R_k$ for $k\geq 26$, then, since $p_j\in S_{25}$, 
$\tilde\Ccal_j$ is also a component of the entire graph $\Dcal_\eta(V,r)$ which has diameter $< 1000r$ and contains $p_j$, contradicting the 
earlier observation that the diameter of the component that contains $p_j$ is $>1000 r$). See Figure~\ref{fig:onecomp} (the left part).
\begin{figure}[!hbt]
\begin{center}
\input{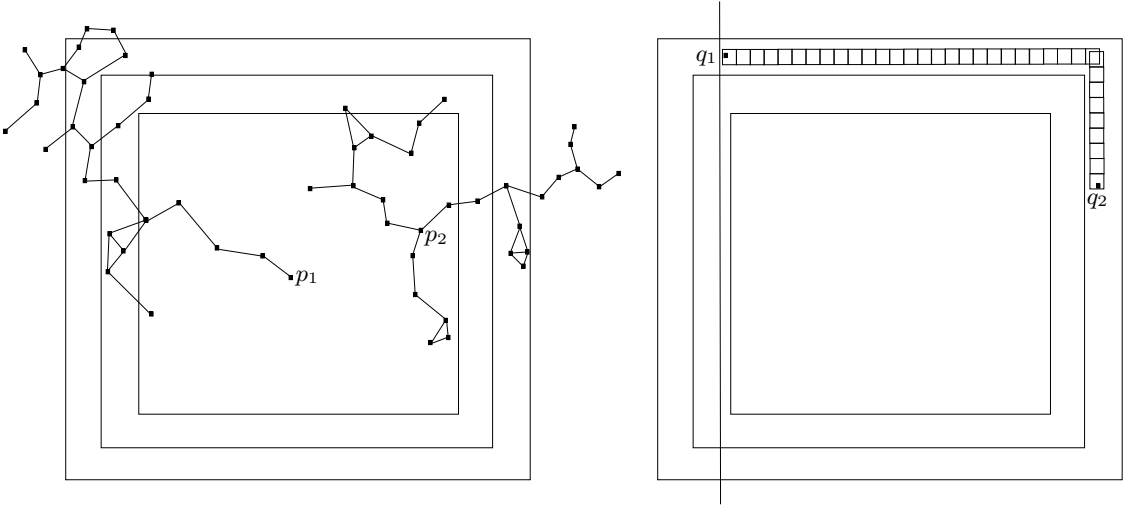}
\end{center}
\caption{Two points $p_1,p_2$ at small distance, in large components, but without a short path between them.}
\label{fig:onecomp}
\end{figure}
Pick an arbitrary $26\leq k\leq 70$.
Let $q_1\in\tilde\Ccal_1, q_2\in\tilde\Ccal_2$ be two points inside $R_k$. 
Provided that $(p_1)_x,(p_1)_y\not\in ((k-1)r,(k-1+1/2\sqrt{2})r)\cup (1-(k-1+1/2\sqrt{2})r, 1-(k-1)r)$, 
it is easy to construct a sequence of squares $T_1, \dots, T_m \subseteq R_k\cap[0,1]^2$, each of side length $r'/2\sqrt{2}$, such that $q_1\in T_1$ and $q_2\in T_m$, and
$T_{i}\cap T_{i+1} \neq \emptyset$ for all $i=1,\dots,m-1$ (see Figure~\ref{fig:onecomp}, the right part).
Observe that every point of $T_{i}$ is at distance $\leq r'$ of every point of $T_{i+1}$.
Hence, if every $T_i$ were to contain at least one point of $\Dcal_\eta(V,r)$, then 
there would be a path between $q_1$ and $q_2$ in $\tilde{\Dcal}$. 
But this cannot be since $\tilde\Ccal_1$ and $\tilde\Ccal_2$ are distinct components of $\tilde{\Dcal}$.
Hence, 
for each $26\leq k\leq 70$ for which $(p_1)_x,(p_1)_y\not\in((k-1)r,(k-1+1/2\sqrt{2})r)\cup (1-(k-1+1/2\sqrt{2})r, 1-(k-1)r)$ 
(note there are at most
$2$ values of $k$ for which this fails), there is at least one
square $T\subseteq R_k$ of side length $r'/2\sqrt{2}$ that does not contain any dense point.

For every $26\leq k\leq 70$ for which this is possible, pick such a square, let $A$ denote the union of these squares, and
set $\Scal := A\cap\Hcal_\eta(r)$.
Clearly $\diam(A) < 10^5r$ and $A$ is a Boolean combination of less than 1000 halfplanes and balls of radius $\leq r'$.
Lemma~\ref{lem:Sarea} and~\eqref{eq:etadef} thus give:
\[ \begin{split} 
|\Scal| 
& \geq 
\area( A ) / (\eta r)^2 - C\eta^{-1} \\
& = 
\frac{43}{8}(1-\eta\sqrt{2})^2\eta^{-2} - C\eta^{-1} \\
& > 
(1+\eps)\pi\eta^{-2}.
\end{split} \]
But then some $q\in\Scal$ must be dense, contradiction.

\noindent {\bf Proof that \pp{6} holds:}
Let us call a point $p\in\Dcal_{\eta}(V,r)$ {\em large} if it is in a component of diameter $\geq r'$.
We first claim that any square $A\subseteq [0,1]^2$ of side length $5r$ contains 
a large point.
To see this, pick such a square $A$, remove a vertical strip of width $r$ from the middle, and denote the
two remaining rectangles of dimensions $2r\times 5r$ by $A_1,A_2$.
Let $\Scal_j := A_j \cap \Hcal_\eta(r)$ for $j=1,2$. 
Then, by Lemma~\ref{lem:Sarea} and~\eqref{eq:etadef} we have
$|\Scal_j| \geq 10\eta^{-2} - C\eta^{-1} > (1+\eps)\pi\eta^{-2}$.
Hence, each $\Scal_j$ contains at least one dense point. A dense point in $\Scal_1$ and a dense point in $\Scal_2$ have distance
between $r$ and $5r\sqrt{2}$, so by \pp{2} at least one of them is large.
So the claim holds.

Now pick two arbitrary large points $p_1,p_2$ of $\Dcal_\eta(V,r)$.
It is easy to construct a sequence $T_1,\dots,T_m$ of squares of side $5r$ such that
$T_i \subseteq [0,1]^2$ for $i=1,\dots,m$, $p_1\in T_1, p_2\in T_m$ and 
$T_i\cap T_{i+1}\neq \emptyset$ for $i=1,\dots,m-1$.
Observe that any point in $T_i$ and any point in $T_{i+1}$ have distance $<25r$. 
By \pp{5} every large point
of $T_i$ is in the same component as every large point in $T_{i+1}$, and, since every $T_i$ has at least one large point, this gives that $p_1$ and $p_2$ lie in the same component.

Since $p_1,p_2$ were arbitrary large points, this shows that all large points lie in the same component.
There is at least one large point (inside any square of side $5r$), so that there indeed is exactly 
one component of diameter $\geq r'$.
\end{prf}  %%% of Proposition prop:Hstruct

It now remains to prove Theorem~\ref{thm:deterministic}. We will make use of the following observation that is essentially to be found in~\cite{diazmitscheperez07}, but is not stated explicilty there. For completeness we include the (short) proof.

\begin{lemma}
Any connected unit disk graph has a spanning tree of maximum degree $\leq 26$.
\label{lem:UDspanningtree}
\end{lemma}

\begin{prf} Let $G=G(V,r)$ be a connected unit disk graph.
For $i,j\in\Zed$ set $V_{i,j} := V \cap [ir/\sqrt{2}, (i+1)r/\sqrt{2})\times [jr/\sqrt{2}, (j+1)r/\sqrt{2})$.
Observe that the vertices of $V_{i,j}$ form a clique in $G$ for each $i,j$, and that 
there can be an edge $vw$ in $G$ between $w\in V_{i,j}$ and $v\in V_{k,l}$ only if $|i-k|,|j-l|\leq 2$.
We construct a subgraph $T$ of $G$ as follows:
\begin{itemize}
\item For each $i,j$ such that $V_{i,j}\neq\emptyset$ we delete all edges between distinct vertices of $V_{i,j}$ except for a path going through all its vertices.
\item For each pair $(i,j)\neq(k,l)$ such that there exists an edge between a vertex in $V_{i,j}$ and a vertex in $V_{k,l}$
we delete all but one of these edges.
\end{itemize}
Observe that if $vw$ is an edge of $G$ then there is a $vw$-path in $T$. Hence $T$ is a spanning 
subgraph of $G$.
It is also clear that $T$ has maximum degree at most $26$, because any vertex is joined to at most 2 vertices in the same $V_{i,j}$ and at most 24 vertices in different $V_{i,j}$s.
If $T$ is not a tree then we can delete additional edges to make it into a tree.
\end{prf}

\begin{prf}[ of Theorem~\ref{thm:deterministic}]
Let $0 < \eta < 1/\sqrt{2}, V \subseteq\eR^2$ and $r>0$ be such that $\Hcal_\eta(r), \Dcal_\eta(V,r)$ and $\Bcal_\eta(r)$ 
satisfy properties \pp{1}-\pp{6} and suppose that $r\leq\rho\leq 2r$ is such that $G(V,\rho)$ is 2-connected.
Let us enumerate the components of $\Dcal_\eta(V,r)$ and the components of $\Bcal_\eta(V,r)$ as $\Ccal_1,\dots, \Ccal_m$, where 
$\Ccal_1$ is the unique component of $\Dcal_\eta(V,r)$ of geometric diameter $\geq 1000r$ and all other
$\Ccal_i$ have diameter $< r'$
(observe that by \pp{4} all components of
$\Bcal_\eta(V,r)$ have geometric diameter $<r'$).
% and that the components of $\Bcal(V,r)$ are vertex-disjoint from those of $\Dcal_\eta(V,r)$ by definition of a bad vertex).
%Let $\Tcal \subseteq \Ccal_1$ be a spanning tree of $\Ccal_1$.
For $p\in\Hcal_\eta(r)$ denote $V_p := V \cap (p+[0,\eta r)^2)$ and 
for $i=1,\dots,m$ let us set $V_{\Ccal_i} := \bigcup_{p\in\Ccal_i} V_p$.
Observe that if $pq$ is an edge of $\Hcal_\eta(r)$ then
$vw$ is an edge of $G(V,\rho)$ for all $v\in V_p, w\in V_q$
(if $v_1\in V_{p_1}, v_2\in V_{p_2}$ with $\norm{p_1-p_2} < r' = r(1-\eta r\sqrt{2})$ then 
$\norm{v_1-v_2} \leq \norm{p_1-p_2} + \norm{(v_1-p_1)-(v_2-p_2)} < r' + \eta r\sqrt{2} = r$).
Amongst other things this shows that $V_{\Ccal_i}$ induces a clique in $G(V,\rho)$ for $i=2,\dots,m$.

\begin{claim}
For each $i=2,\dots,m$ for which $|V_{\Ccal_i}| > 0$ there are paths $P_1^i, P_2^i$ in $G(V,\rho)$ such that:
\begin{enumerate} 
\item $P_1^i, P_2^i$ both have one endvertex in $V_{\Ccal_i}$ and one endvertex in $V_{\Ccal_1}$ and 
all their other vertices in $V\setminus \bigcup_{j=1}^m V_{\Ccal_j}$;
\item $P_1^i$ and $P_2^i$ are vertex-disjoint if $|V_{\Ccal_i}| \geq 2$ and if
$V_{\Ccal_i} = \{v\}$ they share only the vertex $v$ but no other vertices;
\item\label{itm:clean3} There is a $p \in \Ccal_i$ such that both $P_1^i$ and $P_2^i$ are contained in the disk
$B(p,6r)$.
\end{enumerate}
\label{clm:cleanuppaths}
\end{claim}

\begin{prf}[ of Claim~\ref{clm:cleanuppaths}]
If $|V_{\Ccal_i}|\geq 2$ then, since $G(V,\rho)$ is 2-connected, we can pick 
distinct vertices $a_1,a_2\in V_{\Ccal_1}$ and distinct $b_1,b_2\in V_{\Ccal_i}$ and a $a_1b_1$-path $P_1$ and
a $a_2b_2$-path $P_2$ such that $P_1$ and $P_2$ are vertex-disjoint (exercise 4.2.9 on page 173 of~\cite{westboek}).
If $|V_{\Ccal_i}| = \{b_1\}$, then 
we can pick distinct vertices $a_1,a_2\in V_{\Ccal_1}$
and a $a_1b_1$-path $P_1$ and a $a_2b_1$-path $P_2$ 
whose only common vertex is $b_1$ (exercise 4.2.8 on page 173 of~\cite{westboek}).
If $|V_{\Ccal_i}|=1$ then we set $b_2=b_1$ in the rest of the proof.

By switching to subpaths if necessary, we can assume that $a_j$ is the only vertex of $V_{\Ccal_1}$ on $P_j$ 
and $b_j$ is the only vertex of $V_{\Ccal_i}$ on $P_j$ for $j=1,2$.
Let $p_1,p_2\in\Ccal_i$ be such that $b_j\in V_{p_j}$ for $j=1,2$.
We will now show that we can assume that $P_j \subseteq B(p_j; 5r)$ for $j=1,2$
(which implies that both are contained in $B(p_1, 6r)$ as $\norm{p_1-p_2} \leq r' < r$).

Suppose that $P_1$ is not contained in $B(p_1;5r)$.
Write $P_1 = w_0w_1w_2 \dots w_k$ where $w_0=b_1$ and $w_k=a_1$.
Let $j$ be the first index such that $\norm{w_j - p_1} > 2r$.
Observe that 
\[ 
\norm{w_j-p_1} \leq \norm{w_{j-1}-p_1} + \rho %\leq (1+\frac{1}{100}+\eta\sqrt{2})r + \rho 
\leq 4r.
\]
Let $p\in \Hcal_\eta(r)$ be such that $w_j \in V_p$.
Since $\norm{w_j-p_1} - \norm{p-w_j} \leq \norm{p-p_1} \leq \norm{w_j-p_1} + \norm{p-w_j}$
and $\norm{w_j-p} < \eta r\sqrt{2}$ we have
\[ 
2r' < \norm{p-p_1} < 5r.
\]
Depending on whether $\Ccal_i$ is a small component of $\Dcal_\eta(V,r)$ or a component of $\Bcal_\eta(V,r)$, 
by either \pp{3} or \pp{4} we have that $p$ cannot be bad.
Hence there is a dense $q\in\Dcal_\eta(V,r)$ with $\norm{p-q} \leq r'$.
Observe that $\norm{q-p_1} \geq \norm{p-p_1} - 
\norm{q-p} > r'$.
Either \pp{2} or \pp{3} (depending on whether $\Ccal_i$ is a a component 
of $\Bcal_\eta(V,r)$ or a small component of $\Dcal_\eta(V,r)$) now gives that $q\in\Ccal_1$.
Let us pick an $a_1' \in V_q$ that is distinct from
$a_2$ (since $q$ is dense such a $a_1'$ certainly exists).
Then $w_ja_1'$ is an edge of $G(V,\rho)$, and 
\[ \norm{a_1'-p_1} \leq \norm{w_j-p_1}+\norm{a_1'-w_j} \leq 5r.
\]
Hence the path $P_1' = v_1 w_1 \dots w_{j}a_1'$ is as required.
The same argument shows that we can also assume that 
$P_2 \subseteq B(p_2,5r)$.
\end{prf}

\noindent
Part~\ref{itm:clean3} of Claim~\ref{clm:cleanuppaths} implies the following:

\begin{claim}
$P_j^i$ and $P_{j'}^{i'}$ are vertex disjoint for all $i\neq i' \in\{2,\dots,m\}$ and $j,j'\in\{1,2\}$.
\label{clm:Pjidisj}
\end{claim}

\begin{prf}[ of Claim~\ref{clm:Pjidisj}]
Suppose there exists a common vertex $v$. By part~\ref{itm:clean3} of Claim~\ref{clm:cleanuppaths}
there exist $p\in \Ccal_i, p'\in\Ccal_{i'}$ with $\norm{p-p'} \leq \norm{p-v}+\norm{p'-v} < 12r$. 
But this contradicts either \pp{2}, \pp{3} or \pp{4}, depending on what kind of components
$\Ccal_i, \Ccal_{i'}$ are.
\end{prf}

For $i=2,\dots,m$, and $j=1,2$, let $a_j^i$ denote the endpoint of $P_j^i$ in $V_{\Ccal_1}$ and let
$b_j^i$ denote the endpoint of $P_j^i$ in $V_{\Ccal_i}$, and
let $p_1^i, p_2^i \in \Ccal_1$ be such that $a_j^i\in V_{p_j^i}$.
Since $\norm{p_1^i-p_2^i} < 25r$, there is a $p_1^ip_2^i$-path $\Pcal_i$ in $\Ccal_1$ such that 
$\Pcal_i \subseteq B(p_1^i, 100r )$ by \pp{5}.

\begin{claim}
If $i\neq i'$ then $\Pcal_i$ and $\Pcal_{i'}$ are vertex-disjoint.
\label{clm:Pcaldisj}
\end{claim}

\begin{prf}[ of Claim~\ref{clm:Pcaldisj}]
Suppose that some $q\in\Dcal(V,r)$ lies on both paths.
There is a $p\in\Ccal_i$ such that 
$\norm{p-p_1^i} \leq \norm{p-a_1^i}+\norm{a_1^i-p_1^i} \leq 6r +\eta r\sqrt{2}$.
Hence $\norm{p-q} < \norm{q-p_1^i} + 7r \leq 107r$. Similarly, there is a $p'\in\Ccal_{i'}$ such that $\norm{p'-q} < 107r$.
But then $\norm{p-p'} < 214r$, which contradicts one of \pp{2}, \pp{3} or \pp{4} (depending on what kind of components
$\Ccal_i, \Ccal_{i'}$ are).
\end{prf}

To each vertex 
$v\in V\setminus\left(\bigcup_{i=1,\dots,m} V_{\Ccal_i} \cup \bigcup_{i=2,\dots,m,\atop j=1,2}P_j^i\right)$ we will
attach a label as follows.
For such a $v$ there is a $p\in\Hcal_\eta(r)$ such that 
$v\in V_p$. Note that $p$ cannot be bad (otherwise we would have $v\in V_{\Ccal_i}$ for some $i$).
Hence there is at least one dense $q\in\Dcal_\eta(V,r)$ with $\norm{p-q'} < r'$.
Pick an arbitrary such $q$ and label $v$ with $q$ 
(note $vw$ is an edge of $G(V,\rho)$ for all $w\in V_q$).
For a dense $q\in\Dcal_\eta(V,r)$ let us set $L_q := \{ w \in V : w \text{ is labelled } q \}$ and for
$i=1,\dots,m$ will write $L_{\Ccal_i} := \bigcup_{q\in\Ccal_i} L_q$.

Let us observe that for any dense $q\in\Dcal_\eta(V,r)$ and any 7 points $v_1,\dots,v_7 \in V_q$ there is a $v_1v_7$-path
that contains the vertices of $L_q$ and the vertices $v_1,\dots,v_7$ but no 
other vertices. 
This is because all vertices labelled $q$ are adjacent to 
all vertices of $V_q$ and the vertices labelled $q$ can be partitioned into 6 cliques, since
the vertices labelled $q$ all lie inside the disc $B(q,r')$ and this disk can be dissected into 6 sectors of 60 
degrees (each of which has geometric diameter $r'$)-- see Figure~\ref{fig:sixenough}.
We will call such a path a {\em clean-up path} (at $q$).

\begin{figure}
\begin{center}
\includegraphics{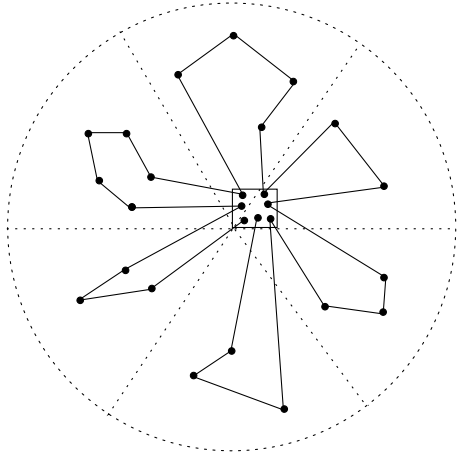}
\end{center}
\caption{A clean-up path}
\label{fig:sixenough}
\end{figure}

\begin{claim}
For $i=2,\dots,m$ and each $q\in\Pcal_i$ and every pair $v,w\in V_q$, there exists
a $vw$-path $P_{v,w}^i$ in $G(V,\rho)$ that visits all vertices of $P_1^i, P_2^i, V_{\Ccal_i}$ and $L_{\Ccal_i}$, 
and at most four vertices from $V_p$ for each $p\in\Pcal_i$, but no other vertices.
\label{clm:Q}
\end{claim}

\begin{prf}[ of Claim~\ref{clm:Q}]
Let us write $\Pcal_i = q_1\dots,q_N$, where $q=q_j$ for some $1\leq j\leq N$ and
$a_1^i\in V_{q_1}, a_2^i\in V_{q_N}$.
We can assume that $v\neq a_2^i$ and $w\neq a_1^i$, by relabelling if necessary.

First suppose that  $b_1^i=b_2^i$. 
In this case we must have $V_{\Ccal_i} = \{b_1^i\}$.
But then $\Ccal_i$ must consist of bad points and $L_{\Ccal_i} = \emptyset$.
We construct the path $P=P_{v,w}^i$ as follows.
Starting from $v$ we go to a vertex $v_{j-1}\in V_{q_{j-1}}$, from there
to $v_{j-2}\in V_{q_{j-2}}$ and so on until $v_1\in V_1$, where we make sure
to pick $v_1=a_1^i$. Next we follow $P_1^i$ to $b_1^i=b_2^i$.  
(If $j=1$ and $v=a_1^i$ then we immediately embark on $P_1^i$. If $j=1$ and $v\neq a_1^i$ then
we first move from $v$ to $a_1^i$ and then embark on $P_1^i$).
Now we follow $P_2^i$ to $a_2^i$. If it happens that $j=N$ and $w=a_2^i$ then we are done. 
If $j=N$ and $w\neq a_2^i$ we jump from $a_2^i$ to $w$ and we are done.
Otherwise we move from $a_2^i$ to a $v_{N-1} \in V_{q_{N-1}}$,
from there to a $v_{N-2}\in V_{q_{N-2}}$ and so on until
$v_j \in V_{q_j}$, where we make sure to pick $v_j=w$. 

Now assume $b_1^i\neq b_2^i$.
For each $q\in\Ccal_i$ such that $|L_q| > 0$, we pick
7 vertices in $V_q$ different from $b_1^i,b_2^i$ (there exist 7 such vertices, because $q$ occurs as a label and is therefore dense) 
and construct the corresponding clean-up path. 
We now construct the path $P_{v,w}^i$ as follows.
We start by going from $v$ to $b_1^i$ in the same way as above.
Since $V_{\Ccal_i}$ is a clique, we can start from $b_1^i$, jump to an endvertex of the first clean-up path, follow it,
jump from its other endvertex to an endvertex of another clean-up path, follow that path and so on until the last clean-up path.
We then follow a path trough the remaining vertices of $V_{\Ccal_i}$, arriving at $b_2^i$. Finally
we follow $P_2^i$ to $a_2^i$, and go from $a_2^i$ back to $w$ in the same way as above.
\end{prf}

By Lemma~\ref{lem:UDspanningtree} there exists a spanning tree $\Tcal$ of $\Ccal_1$ with maximum degree 
at most 26.
Let $\Wcal = q_0\dots q_N$ (with $q_0=q_N$) be a closed walk on $\Tcal$ that traverses every edge exactly twice (once in each direction).
Such a walk can for instance be obtained by tracing the steps of a depth-first search algorithm on $\Tcal$.
Observe that $\Wcal$ visits each node $q\in\Ccal_1$ at most 26 times, since the maximum degree of $\Tcal$ is at most 26.
We shall now describe a construction of a Hamilton cycle in $G(V,\rho)$.
It is convenient to consider "timesteps" $t=0,\dots,N$, where we envisage
ourselves performing the walk $\Wcal$ while at the same time constructing the cycle $C$.
At the beginning of timestep $t$, the cycle $C$ under construction is always at a vertex $v\in V_{q_t}$
and at the end of timestep $t<N$ we are at a vertex $w\in V_{q_{t+1}}$. 
We start the cycle from an arbitrary vertex $\alpha_0\in V_{q_0}$.
At the beginning of timestep $t$ we are in some vertex $v\in V_{q_t}$. We now apply the following rules at each timestep $t=0,\dots,N$:
\begin{itemize}
\item[\rr{1}] If it is the first time $\Wcal$ visits $q_t$
(i.e.~$q_t$ is distinct from $q_0, \dots, q_{t-1}$), 
and $q_t$ lies on $\Pcal_i$ for some $i=2,\dots,m$, and 
it is the first vertex of $\Pcal_i$ that occurs on $\Wcal$ then we pick an arbitrary
$w\in V_{q_t}\setminus\{v\}$ and we follow the path $P_{v,w}^i$ from $v$ to $w$.
This timestep has not finished yet. We next apply either Rule 2 or Rule 3 (whichever applies).

\item[\rr{2}] If it is not the last time $\Wcal$ visits $q_t$, then we simply pick a not yet visited $w\in V_{q_{t+1}}$ and
go there. {\bf End of timestep $t$}.

\item[\rr{3}] If it is the last time that $\Wcal$ visits $q_t$
(i.e.~$q_t$ is distinct from $q_{t+1},\dots,q_N$)
then we do the following.
We are currently in a vertex $v\in V_{q_t}$.
Pick vertices $v_1,\dots,v_6 \in V_{q_t}$ that have not been visited yet, and
follow a clean-up path between $v$ and $v_6$ that visits all vertices labelled $q_t$ and
$v_1,\dots,v_6$, but no other vertices. 
Now we continue by visiting all vertices of $V_{q_t}$
that have not yet been visited. Finally, provided $t<N$, we pick a not yet visited $w\in V_{q_{t+1}}$
and go to $w$. If $t=N$, then we go to the initial vertex $\alpha_0$, completing the cycle.
{\bf End of timestep $t$}.
\end{itemize}
Let us now explain why this construction works.
At each timestep $t$ the Rules 1-3 require unused vertices in $V_{q_t}$, so we need to argue amongst other things 
that we never run out of vertices.
Pick an arbitrary $q\in\Ccal_1$.
By Claim~\ref{clm:Pcaldisj} there is at most one $2\leq i\leq m$ such that $q\in\Pcal_i$.
Rule 1 is applied exactly once to a $p\in\Pcal_i$, and when that happens at most 
4 new vertices of $V_q$ are used.
Rule 2 is applied at most 25 times to $q$, and each time one new vertex of $V_q$ is used.
Thus, at the start of the timestep when $\Wcal$ visits $q$ for the last time, at least $100-4-25 = 71$ vertices of $V_q$ are left, which is
more than enough to construct the clean-up path. 
So we never get stuck.

%It is also clear that $C$ is a cycle (we never visit a vertex of $v\in V$ twice, ).
We still need to argue that our construction produces a Hamilton cycle.
Recall that $V$ can be partitioned into the sets $V_{\Ccal_i}, L_{\Ccal_i} : i=1,\dots,m$ and  
$P^i_j \setminus \{a_j^i, b_j^i\} : i=2,\dots,m, j=1,2$ (by construction of
$V_{\Ccal_i}$ and $L_{\Ccal_i}$ and by Claim~\ref{clm:Pjidisj}).
Consider an arbitrary $v\in V$.
If $v \in P^i_j\setminus \{a_j^i, b_j^i\}$ for some $i=2,\dots,m, j=1,2$,
then $C$ visits $v$ exactly one, namely in the time step when 
$\Wcal$ first visits a vertex of $\Pcal_i$.
Similarly, if $v\in V_{\Ccal_i}$ or if $v\in L_{\Ccal_i}$ for some $i=2,\dots,m$, then 
$C$ visits $v$ exactly once, namely in the time step when 
$\Wcal$ first visits a vertex of $\Pcal_i$.
If $v\in L_{\Ccal_1}$ then $C$ visits it exactly once, namely 
at the timestep when $\Wcal$ visits $q$ for the last time where $q\in\Ccal_1$ is such that $v$ is labelled $q$.
It is also clear that $C$ visits every vertex $v\in V_{\Ccal_1}$ exactly once (in Rules 1 and 2
we always take new vertices from $V_q$, and when Rule 3 is finally applied to $q$ we make sure to visit 
all remaining vertices of $V_q$).
Thus, $C$ visits every $v\in V$ exactly once and, since in the very end we 
reconnect to the initial vertex $\alpha_0$, it is a Hamilton cycle as required.
\end{prf}

\section{Extension to other norms and higher dimensions}
\label{sec:generalise}

In this section we shall briefly sketch the changes needed to make proof of Theorem~\ref{thm:main} work 
in the case when $X_1,X_2,\dots$ are independent, uniform random points from $[0,1]^d$ with $d\geq 2$ arbitrary 
and when $\norm{.}$ in the definition of the random geometric graph is the $l_p$-norm for some $1<p\leq \infty$.
That is,

\begin{theorem}
For any $d\geq 2$ and $1<p\leq \infty$ the following holds.
If we pick $X_1,\dots, X_n \in [0,1]^d$ i.i.d.~uniformly at random and we 
add the edges $X_iX_j$ by order of increasing $l_p$ norm of $X_i-X_j$ then, with
probability tending to 1 as $n\to\infty$, the resulting graph gets its first
Hamilton cycle at precisely the same time it loses its last vertex of degree
less than two. 
\end{theorem}

We should perhaps remark that the restriction to the $l_p$-norm with $1<p\leq \infty$ is needed 
only because it is imposed by the results of Penrose that we invoke in our proofs (cf.~Theorem 8.4 and 13.17 of~\cite{penroseboek}).
These results of Penrose show a notable difference between the case when the points $X_1,\dots, X_n$ are chosen uniformly at random 
from the unit hypercube and the case when they are chosen from the $d$-dimensional torus (i.e.~if we identify opposite facets
of the unit hypercube). The restriction to $l_p$-norms with $1<p\leq \infty$ is imposed
by Penrose for the unit hypercube (but not for the torus) to deal with the technical difficulties that arise from 
``boundary effects''.

Most of the proofs go through almost unaltered if we change the relevant constants etc. in the following way.
When a square appears  in the proofs for the 2-dimensional, Euclidean case, it should usually be replaced by a $d$-th power. 
Instead of $\area(.)$ we need to put $\vol(.)$, the $d$-dimensional volume.
Whenever the constant $\pi$ occurs it should be replaced by $\theta := \vol( B(0,1) )$, the volume of 
the unit ball wrt.~the $l_p$-norm.
Wherever the constant $\sqrt{2}$ appears, it should be replaced by $d^{1/p} = \diam([0,1]^d)$, the diameter 
of the $d$-dimensional hypercube as measured by the $l_p$-norm (here we interpret $1/\infty$ as $0$, so that $d^{1/\infty} = 1$). 
For example, we now set $r' := r(1-\eta d^{1/p})$.
Instead of the numbers $10^5, 1000, 100, 25$ we put suitably chosen large constants.
In particular, in the definition of $\Hcal_\eta(r), \Dcal_\eta(V,r), \Bcal_\eta(V,r)$ a point of
$p \in \Hcal_\eta(r)$ is dense if the cube $p+[0,\eta r)^d$ contains at least $K$ points for a suitably chosen 
constant $K$ (that will have to be larger than 100 for some choices of $d,p$).

In the statement and proof of (the analogues of) Proposition~\ref{prop:Hstruct} and Lemma~\ref{lem:Sdense} for the general case 
we can put 
\[ r_n := \left( (1-\delta)\left(\frac{2^{d-1}}{d\theta}\right)\ln n / n\right)^{1/d}, \]
where $\delta = \delta(\eps)$ is a suitably chosen small constant.
It can be read off from Theorem 8.4 in~\cite{penroseboek}, together with Theorem 13.17 in~\cite{penroseboek} 
(the version of Theorem~\ref{thm:penrosekconn} for arbitrary dimension and the $l_p$-norm), 
that $r_n < \rho_n(\text{2-connected}) < 2r_n$.
We should perhaps remark that, although it is possible to have $\delta$ tend to 0 in a suitable way, 
it cannot be disposed of altogether.
This is because (c.f.~Theorem 8.4 in~\cite{penroseboek}) the
last vertex of degree $< 2$ disappears when $r = \left( (\frac{2^{d-1}}{d\theta}\ln n + c_{d,p}\ln\ln n + O(1)) / n \right)^{\frac{1}{d}}$
where $c_{d,p}$ is a constant that is negative for some choices of $d,p$.

In the higher-dimensional analogue of Lemma~\ref{lem:Sdense} we need to 
distinguish additional subcases to deal with the situation when $\Scal$ is close
to a $k$-dimensional face of $[0,1]^d$, for $k=1,\dots,d-1$.
Let $\sde_k(s)$ denote the set of all $z\in[0,1]^d$ that have $k$ coordinates in $[0,s)\cup(1-s,1]$.
Then $|\Hcal_\eta(r_n) \cap \sde_k( Kr_n)| = O( r_n^{-(d-k)} )$, and also
$|\Ucal_k| = O( r_n^{-(d-k)} ) = O( (n/\ln n)^{(d-k)/d} )$, where $\Ucal_k$ is the collection of all 
sets $\Scal \subseteq \Hcal_\eta(r)\cap\sde_k(Kr)$ with diameter at most $Kr$.
The argument in the proof of Lemma~\ref{lem:Sdense} thus shows that each $\Scal\in\Ucal_k$
with $|\Scal| > (1+\eps)\frac{d-k}{2^{d-1}}\theta\eta^{-d}$ contains a dense point.

Lemma~\ref{lem:Sarea} and its proof essentially go through unaltered if we replace $\area(.)$ by the $d$-dimensional volume $\vol(.)$, $\sqrt{2}$ by $d^{1/p}$, $\eta^{-2}$ by $\eta^{-d}$ and
$\eta^{-1}$ by $\eta^{-(d-1)}$. 

In the proof of \pp{1}, we now 
pick $2d$ vectors from $\Ccal$, two for each coordinate. 
For $i=1,\dots,d$ we let $p_i^{-}$ resp.~$p_i^{+}$ be a point of smallest resp.~largest
$i$-th coordinate.
We set $A := \bigcup_i B_{i}^{-}(p_{i}^{-},r')\cup B_{i}^{+}(p_{i}^{+},r')$, where
$B_{i}^{-}(z,s) := \{ z' \in B(z,s) : (z')_i < z_i \}, B_{i}^{+}(z,s) := \{ z' \in B(z,s) : (z')_i > z_i \}$,
and we set $\Scal := A \cap \Hcal_\eta(r_n)$.
Again it is clear that $\Scal$ cannot contain any dense point if $\Ccal$ is a component.  
This time here must exist an $1\leq j \leq d$ such that $p_{j}^{+}-p_{j}^{-} > r' / d^{1/p}$.
We now need to consider the case when one of these $2d$ points is in 
$\sde_k( r' )$ but none lies in $\sde_{k+1}( r')$.
Wlog.~suppose the points are close to the face $\{ z\in[0,1]^2 : z_1=\dots=z_k=0\}$.
First suppose that $j \leq k$. The we can assume (w.l.o.g.) that $j=1$.
We see that $A\cap[0,1]^2$ contains $A_1$ and $A_2$, where
\[ 
\begin{array}{l}
A_1 := \{ z \in B(p_{1}^{+}, r' ) : z_i > (p_{1}^{+})_i \text{ for } i=1,\dots, k\}, \quad \text{ and } \\
A_2 := \{ z \in B(p_{2}^{+}, r' ) : (p_{1}^{-})_1 < z_1 < (p_{1}^{+})_1, \text{ and }
z_i > (p_{2}^{+})_i \text{ for } i=2,\dots, k\}.
\end{array} \]
Note that $A_1$ and $A_2$ are disjoint, that $\vol(A_1) = \theta (r')^d / 2^k$ and that $\vol(A_2) > \theta (r')^d / 2^kd^{1/p}$.
Hence 
\[ \begin{split}
|\Scal| 
& \geq 
\vol(A\cap[0,1]^d) / (\eta r_n)^d - C \eta^{-(d-1)} \\
& \geq  
(1+1/d^{1/p})(1-\eta d^{1/p})^d \theta\eta^{-d} / 2^k - C \eta^{-(d-1)} \\
& > (1+\eps)\theta\eta^{-d} / 2^{-k} \\
& \geq (1+\eps)\frac{d-k}{2^{d-1}}\theta\eta^{-d}, 
\end{split}
\]
(provided $\eps,\eta$ were chosen appropriately) so that $\Scal$ must contain a dense point.

Now consider the case when $j > k$. Then $A\cap [0,1]^2$ contains
\[ \begin{array}{l}
A_1 := \{ z \in B_{j}^{+}(p_{j}^{+}, r' ) : z_i > (p_{j}^{+})_i \text{ for } i=1,\dots, k\}, \\
A_2 := \{ z \in B_{j}^{-}(p_{j}^{-}, r' ) : z_i > (p_{j}^{+})_i \text{ for } i=1,\dots, k\},  \quad \text{ and } \\
A_3 := \{ z \in B(p_{1}^{+}, r' ) : (p_{j}^{-})_{j} < z_{j} < (p_{j}^{+})_j, \text{ and }
z_i > (p_{1}^{+})_i \text{ for } i=1,\dots, k\}.
\end{array} 
\]
Now $A_1, A_2$ both have volume $\theta (r')^d / 2^{k+1}$ and $A_3$ has volume 
at least $\theta (r')^d / 2^k d^{1/p}$. So again $\Scal$ must contain a dense point.

The arguments that reduce \pp{2}-\pp{4} to the proof of \pp{1} work in the same way for other dimensions and norms.
In the proof of \pp{5} we now show that (with $C_1<C_2$ suitable constants) if two points $p_1,p_2$ have distance less than $C_1$ and 
there is no path between them that stays inside $p_1 + [-C_2r,C_2r]^d$, then in all but $d$ of the sets 
$p_1 + [-kr,kr]^d \setminus [-(k-1)r, (k-1)r]^d : k=C_1+1,\dots,C_2$ there is a cube of
side $1 / 2d^{1/p}$ without a dense point inside it.  

In the proof of \pp{6} we merely need to replace squares of side $5r$ with hypercubes of side $Kr$ for some suitable constant $K$.

Lemma~\ref{lem:UDspanningtree} and its proof generalise to give that
any connected, (non-random) geometric graph in dimension $d$ with the $l_p$-norm 
has a spanning tree of maximum degree at most $(2\lceil d^{1/p}\rceil + 1)^d +1$.

The proof of Theorem~\ref{thm:deterministic} also generalises with only minor modifications, the most important
one being in the definition of a clean-up path. For any dimension $d$ and $1<p\leq \infty$ there 
exists a finite $k$ such the unit ball wrt.~the $l_p$ norm can be partitioned into
$k$ parts each of diameter $\leq 1$ (covering the ball by hypercubes of side $1/d^{1/p}$ shows for instance 
that we can take $k = (2\lceil d^{1/p}\rceil)^d$.)
We can thus construct clean-up paths at each $q\in\Dcal_\eta(V,r)$
that use $k+1$ vertices from $V_q$.

\section{Concluding remarks}
\label{sec:concl}

In this paper we have shown that, with high probability, the least $r$ for which the random geometric graph $G(n,r)$ is 
Hamiltonian coincides with the least $r$ for which it has minimum degree at least 2.
Recall that a graph is {\em pancyclic} if it has cycles
of all lenghts $3\leq k \leq n$.
As shown by {\L}uczak~\cite{luczak91}, the usual random graph becomes pancyclic at 
exactly the same time it loses its last vertex of degree $<2$.
It is natural to ask whether a similar statement can be shown for the random geometric graph.
As it happens, the answer is yes.
Our proof of Theorem~\ref{thm:main} can be adapted to show:
\begin{theorem}\label{thm:pancyc}
$\Pee\left[ \rho_n( \text{pancyclic} ) = \rho_n(\text{minimum degree} \geq 2 )\right] \to 1$ as
$n\to\infty$.
\end{theorem}
This also implies that in Corollary~\ref{cor:maincor} we can replace the word "Hamiltonian" with "pancyclic".
Let us briefly explain how to adapt our proof of Theorem~\ref{thm:main} to give Theorem~\ref{thm:pancyc}.
First note that the proof would still have gone through if we had defined
$p\in\Hcal_\eta(r)$ to be dense if the corresponding square contains at least 1000 points instead of 100.
We will reconsider the way we constructed the Hamilton cycle $C$, and show that 
for every $1\leq k \leq n-3$ there are $k$ points from $V$ that we can omit and construct a cycle through 
the remaining points applying Rules 1-3 in the same way as in the proof of Theorem~\ref{thm:deterministic}.
For any set of vertices $A\subseteq\bigcup_{i=1}^m L_{\Ccal_i}$
we can construct a cycle through $V\setminus A$ by simply 
omitting the vertices of $A$ and proceeding as in the proof of Theorem~\ref{thm:deterministic}
(the vertices of $A$ will simply be omitted from the corresponding clean up paths).
Thus we have cycles of lengths $n-k$ for $k=0,\dots, \sum_i |L_{\Ccal_i}|$.
Let us now omit all vertices in $\bigcup_{i\geq 1} L_{\Ccal_i}$, and consider $\Ccal_2$.
Since $\Ccal_2$ is a clique, we can omit all vertices except $b_1^2, b_2^2$ one by one and
each time construct a cycle through the remaining points.
Let us omit $V_{\Ccal_2} \setminus \{b_1^2, b_2^2\}$ as well as $\bigcup_{i\geq 1} L_{\Ccal_i}$ in the sequel. 
We can assume w.l.o.g.~that $P_1^2, P_2^2$ have no shortcuts (in other words we can assume they are induced paths), because
we could have easily insisted on this in the proof of Claim~\ref{clm:cleanuppaths}.
This means that each of these paths is the union of two stable sets.
Notice that if $\{v_1,\dots,v_l\} \subseteq B(p,6r)$ is a stable set then the
discs $B(v_1, r/2), \dots, B(v_l, r/2)$ are disjoint and contained in $B(p,6\frac12)$, so that 
$l \leq \pi (6\frac12)^2 / \pi (\frac12)^2 = 169$.
This shows that $P_1^2, P_2^2$ each have at most $338$ vertices.
For each $k=1,\dots,|P_1^2|+|P_2^2|-1$ we can omit $k$ points from $V_p$ for some $p\in\Pcal_2$
(since $|V_p| \geq 1000 > 2\cdot 338$ this can be done). There is always a cycle through the remaining vertices.
Now put back those points from $V_p$ and remove $P_1^2$ and $P_2^2$. Again there is a cycle through all points that 
have not been removed.
Having removed $V_{\Ccal_2}$ and $P_1^2, P_2^2$, we can repeat the same procedure for $i=3,\dots, m$.
We see that there is a
cycle of length $n-k$ for $k=0,\dots, \sum_{i=1}^m |L_{\Ccal_i}|
+ \sum_{i=2}^m |V_{\Ccal_i}| + \sum_{i=2}^m (|P_1^i|+|P_2^i|)$.
Removing $\bigcup_{i=1}^m L_{\Ccal_i}$, $\bigcup_{i=2}^m V_{\Ccal_i}$ and $\bigcup_{i=2}^m P_1^i \cup P_2^i$, we
are only left with vertices of $V_{\Ccal_1}$, and for each $p\in\Ccal_1$ we still have at least 
$1000-2 = 998$ points of $V_p$ left over (it contains at most 2 endpoints of $P_j^i$s).
We omit the remaining points one by one, starting with points in squares corresponding to leafs of $\Tcal$.
Once we have run out of those, we continue with points in squares corresponding to leafs of the subtree 
of $\Tcal$ induced by the nonempty squares, and so on.
We see that are indeed able to construct cycles of all lengths.

Observe that if $\delta(G)$ denotes the minimum degree of the graph $G$, then there can be at most
$\lfloor\delta(G)/2\rfloor$ edge disjoint Hamilton cycles in $G$. 
Bollob\'as and Frieze~\cite{bollobasfrieze85} have shown that, with high probability, the ordinary random graph 
has $k$ edge-disjoint Hamilton cycles for the first time at precisely the same moment
it first achieves minimum degree $2k$.
Perhaps methods similar to ours will prove:
\begin{conjecture}
$\rho_n(\text{there exist } k \text{ edge disjoint Hamilton cycles })
= \rho_n( \text{minimum degree} \geq 2k)$ w.h.p., for any fixed $k \in\eN$.
\end{conjecture}

Let $H_\delta$ denote the graph property that there are $\lfloor \delta(G)/2\rfloor$ edge disjoint Hamilton cycles in the graph $G$.
It has been conjectured (see e.g.~\cite{friezekrivelev05}) 
that $H_\delta$ holds w.h.p.~for all choices of the sequence $(m_n)_n$ in the $G(n,m_n)$ model. 
This is known to be true for choices of $(m_n)_n$ for which $G(n,m_n)$ has 
minimum degree $o(\ln n)$ w.h.p.~(c.f.~\cite{friezekrivelev08}), but it is still open in general.
A natural question is therefore:
 
\begin{question}
Does $H_\delta$ hold w.h.p.~for the random geometric graph $G(n,r_n)$ for all choices of the sequence $(r_n)_n$?
\end{question}

\bibliographystyle{plain}
\bibliography{References}

\end{document}